\documentclass[12pt]{amsart}
\usepackage{amsfonts,amssymb}
\usepackage{epsfig}
\usepackage{latexsym}
\usepackage[all]{xy}
\usepackage{amsmath}
\usepackage{amsfonts}
\usepackage{amssymb}

\def\A{\mathcal A}

\def\O{\mathcal O}
\newtheorem{Def}{Definition}[section]

\newtheorem{Remark}{Remark}[section]

\begin{document}
\title{ An Inverse system of nonempty objects with empty limit}
\author{Satya Deo and Veerendra Vikram Awasthi}
\date{}


\maketitle

\begin{abstract}
In this article we give an explicit example of an inverse system with nonempty sets and onto bonding maps such that its inverse limit is empty.
\end{abstract}

\section{Introduction}
 It has been often quoted, without giving an example, that
there are inverse systems in which all the objects are
nonempty  and the bonding maps are onto, but the inverse limit of the inverse system is
empty. (see e.g.\cite{es}, Dugundji \cite{dug} pg 427, last para). It is hard to believe that there should be an example like this, but the fact is that there exists such an example. The original
paper by L. Henkin \cite{hen} dealing with this problem
gives a theorem  which implies that there are several examples, but
the explanation and proof of the theorem is too abstract to have a
clear idea of a specific example.

\bigskip In this note, we present a concrete example of an inverse
system of sets and maps having the stated property, 
which can be easily understood. It will then easily follow that we can
have  an example of such an inverse system in any category admitting arbitrary products, e.g., in the category of modules
and homomorphisms, category of topological spaces and continuous maps, etc.

\section{Notations and Preliminaries}

\noindent First, let us have the following well known definitions (see \cite{dug}):

{\Def \rm  A binary relation $R$ in a set $A$ is called a {\bf preorder} if it is reflexive and transitive. A set together with a definite preorder is called  a {\bf preordered set}.}

{\Def \rm  Let $\A$ be a preordered set and $\{X_{\alpha}\; |\; \alpha\in \A\}$ be a family of topological spaces indexed by $\A$. For each pair of indices $\alpha, \beta$ satisfying $\alpha < \beta$, assume that there is a given a continuous maps $f^{\beta}_{\alpha}:X_{\beta}\to X_{\alpha}$ and that these maps satisfy the following condition: If $\alpha<\beta<\gamma,$ then $f^{\gamma}_{\alpha}=f^{\beta}_{\alpha}o f^{\gamma}_{\beta}$. Then the family $\{X_{\alpha};\;f^{\beta}_{\alpha}\}$ is called an {\bf inverse system} over $\A$ with topological spaces $X_{\alpha}$ and bonding continuous maps $f^{\beta}_{\alpha}$.}

{\Def \rm  Let $X=\{X_{\alpha,}\;f^{\beta}_{\alpha};\;\alpha\leq\beta\}_{\alpha,
\beta\in\A}$ be an inverse system of topological spaces and
  continuous map $f^{\beta}_{\alpha}:X_{\beta}\to
  X_{\alpha},\;\;\alpha\leq\beta$ based on indexing set $\A$. Consider the product space
  $\prod_{\alpha\in\A}X_{\alpha},$ and let $p_{\alpha}: \Pi
  X_{\alpha}\to X_{\alpha}$ denotes the projection map. Define
$$X'=\{x\in \prod_{\alpha\in\A}X_{\alpha}\;|\;\mbox{\rm
    whenever}\;\alpha\leq\beta,\;f^{\beta}_{\alpha}p_{\beta}(x)=
    p_{\alpha}(x)\}.$$ Then the set $X'$ with subspace topology is
    called the {\bf inverse limit} of the inverse system $X$.} 

{\Remark \rm We have defined  above an inverse system and inverse limit in the category of topological spaces. Clearly, this definition of an inverse system can be made in any category. However, the inverse limit of an inverse system  will exist only in those categories which admit arbitrary products.}

\bigskip\noindent {\bf Ordinal Numbers : }
For definition and well known special properties of ordinal numbers we refer to Dugundji \cite{dug}. 

\bigskip The successor $x^+$ of a set $x$ is defined as $x\cup \{x\}$,
and then $\omega$ was constructed as the smallest set that contains 0
and that contains $x^+$ whenever it contains $x$. Now, the question
arises that what happens if we start with $\omega$, form its successor
$\omega^+$, then form the successor of that, and proceed so on. In
other words, is there something out beyond $\omega, \omega^+,
(\omega^+)^+, \cdots, $ etc., in the same sense in which $\omega$  is
beyond $0,1,2,\cdots,$etc.?

\medskip We mention the
names of some of the first few of them. After $0,1,2,\cdots$ comes
$\omega,$ and after $\omega,\omega +1,\omega +2, \cdots$ comes
$2\omega.$ After $2\omega +1$ (that is, the successor of $2\omega$)
comes $2\omega +2,$ and then $2\omega +3$ ; next after all the terms of the
sequence so begun comes $3 \omega$. At this point another application
of axiom of substitution is required. Next comes $3\omega+1,
3\omega +2,3\omega +3,\cdots ,$ and after them comes $4 \omega$. In
this way we get successively $\omega, 2\omega,
3 \omega,4 \omega,\cdots .$. An application of axiom of substitution
yields something that follows them all in the same sense in which
$\omega$ follows the natural numbers: that something is $\omega^2.$
After the whole thing starts over again : $\omega^2 +1, \omega^2 +2,
\omega^2 +3, \omega^2 +\omega, \omega^2 +\omega+1, \omega^2
+\omega+2,\cdots , \omega^2 +2\omega, \omega^2 +2\omega+1,\omega^2
+2\omega+2, \cdots ,\omega^2 +3\omega, \cdots ,\omega^2 +4\omega,
2\omega^2,\cdots, 3\omega^2, \cdots, \omega^3, \cdots,
\omega^4,\cdots, \omega^{\omega},\cdots,
\omega^{({\omega}^{\omega})},\cdots. $ 

\medskip Since the countable union of countable sets is again countable, each
of the above numbers is countable. Therefore, using the well-ordered
property of  ordinals there exists a smallest ordinal number
$\omega_1$ which contains all of the above numbers and is itself
uncountable. We call $\omega_1$ as the first uncountable ordinal number.

\section{The Example}

 Let $\{0,1,2,\cdots, \omega , \omega +1, \cdots ,
2\omega , 2\omega +1, \cdots, \omega^{2},\omega^{2}+1, \cdots \} $ be
the set of ordinal numbers and $\omega_1$ be the first uncountable
ordinal. Consider the set $\Omega=[0, \omega_1)$ of all ordinals less
  than $\omega_1$. We will construct an
  inverse mapping system $\{X_{\alpha}, f^{\beta}_{\alpha};
  \alpha\leq\beta \}$ based on the directed set $\Omega$ in which all
  the sets $X_{\alpha}$ will be nonempty, the bonding maps
  $f^{\beta}_{\alpha}:X_{\beta}\to X_{\alpha},\;\alpha\leq \beta,$ will be
  onto, but $\varprojlim X_{\alpha}=\phi$.

\bigskip Let us define a {\bf point} to mean a finite sequence of an
even number of 
elements from $\Omega,$ e.g., $$x=(\alpha_1,\alpha_2,\cdots,
\alpha_{2n-1}, \alpha_{2n}),$$ 

which satisfy the following three conditions:   

\begin{enumerate}
\item[(i)] $\alpha_1 < \alpha_2$
\item[(ii)] $\alpha_{2i-1} < \alpha_{2i+2}$ for $0<i<n,$ 
i.e., $\alpha_1<\alpha_4,\;\alpha_3<\alpha_6,\;\alpha_5<\alpha_8,\cdots$
\item[(iii)] $\alpha_{2i+1} < \alpha_{2i+2}$ and $\alpha_{2i+1} \nless
  \alpha_{2j+1}$   for   $0\leq j< i<n,$
\end{enumerate}

\medskip \noindent where $\alpha\nless\beta$ holds when
  neither   $\alpha<\beta$  nor $\alpha=\beta.$ This means \linebreak
  $\alpha_1<\alpha_2,\;\alpha_3<\alpha_4,\;\alpha_5<\alpha_6,\;\alpha_7
  <\alpha_8,\cdots$    and 
 $\alpha_3\nless\alpha_1,\linebreak
  \alpha_5\nless\alpha_3,\;\alpha_5\nless\alpha_1,   \cdots$.

\bigskip
\noindent We may observe that the above conditions  imply that\\ 
$\alpha_1<\alpha_2,\;\alpha_1<\alpha_4,$\\
$\alpha_3<\alpha_4,\;\alpha_3<\alpha_6,$\\
$\alpha_5<\alpha_6,\;\alpha_5<\alpha_8,$\\
$\alpha_7<\alpha_8,\;\alpha_7<\alpha_{10},$ and so on.

\bigskip \noindent We define {\bf index} of the point $x$ given above to be
$\alpha_{2n-1,}$ {\bf order} of $x$ to be $\alpha_{2n}$ and {\bf
  length} of $x$ to be $n$. 

\bigskip \noindent Let us illustrate a few begining sets
$X_{\alpha},\;\alpha\in\Omega$ 

\begin{enumerate}
\item[(i)] $x\in X_0,$ means $x=(0, \alpha)$ where
  $\alpha >0.$ 
\item[(ii)] $x\in X_1,$ means the point $x$ can be one
  of the following   types\\
$x=(1, \alpha)$ where $\alpha >1,$ or\\
$x=(0, \alpha, 1, \beta )$ where $\alpha >0$ and $\beta >1$.
\item[(iii)] $x\in X_2,$ means the point $x$ can be one
  of the following   types\\
$x=(2, \alpha)$ where $\alpha >2,$ or\\
$x=(0, \alpha, 2, \beta )$ where $\alpha >0$ and $\beta >2$ or\\
$x=(1, \alpha, 2, \beta )$ where $\alpha >1$ and $\beta >2$ or\\
$x=(0, \alpha, 1, \beta, 2, \gamma )$ where $\alpha >0$ and $\beta >1$
  and $\gamma >2$ 
\end{enumerate}

\noindent Note that all the points defined above have index 2.
Similarly we can check the elements of  $X_3$ as follows:

\begin{enumerate}
\item[(iv)] $x\in X_3$ means $x$ is one of the following types\\
$x=(3, \alpha)$ where $\alpha >3$ or\\
$x=(0, \alpha, 3, \beta)$ where $ \alpha >0$ and $\beta >3$ or \\
$x=(1, \alpha, 3, \beta)$  where $ \alpha >1$ and $\beta >3$ or \\
$x=(2, \alpha, 3, \beta)\;\mbox{where } \alpha >2$ and $\beta >3$ or \\
$x=(0, \alpha, 1, \beta, 3, \gamma)$ where $\alpha >0$ and
  $\beta >1$ and $\gamma >3$ or \\
$x=(0, \alpha, 2, \beta, 3, \gamma)$ where $ \alpha >0$ and
  $\beta >2$ and $\gamma >3$ or \\
$x=(1, \alpha, 2, \beta, 3, \gamma)$ where $ \alpha >1$ and
$\beta >2$ and $\gamma >3$ or \\ 
$x=(0, \alpha, 1, \beta, 2, \gamma, 3, \delta)$ {where }$ \alpha
>0$ and $\beta >1$ 
and $\gamma >2$ and $\delta>3,\;\;\;\alpha, \beta, \gamma, \delta \in
\Omega.$
\end{enumerate}

\medskip \noindent Thus we have a family of nonempty disjoint sets
$X_{\alpha},\;\;\alpha\in\Omega,$ whose elements are the points with
index $\alpha.$ 

\medskip We now define the bonding maps
$f^{\beta}_{\alpha}:X_{\beta}\to X_{\alpha},\;\;\alpha\leq\beta,
\;\;\alpha, \beta \in\Omega.$  Let $x=(\alpha_1,\alpha_2,\cdots, \alpha_{2n-1},
\alpha_{2n})$ be an arbitrary point in $X_{\beta}$ (so that
$\alpha_{2n-1}=\beta$). We define the image of $x$ in $X_{\alpha}$
under $f^{\beta}_{\alpha}$ as follows:

\medskip\noindent There are two cases :

\bigskip\noindent {\bf Case I:} If $\alpha\leq\alpha_1$, then we define
$f_{\alpha}^{\beta}(x)=(\alpha, \alpha_2)$ and since $x$ is a point in
$X_{\beta},\;\; \alpha_1<\alpha_2$ by condition (i) which implies
$\alpha<\alpha_2$. 
Therefore $(\alpha, \alpha_2)$ is a point with index $\alpha$. Hence
$(\alpha, \alpha_2)\in X_{\alpha}.$

\bigskip\noindent {\bf Case II:} If $\alpha\nless\alpha_1$ then there exist
a least $j,\;0<j<2n-1$ such that $\alpha<\alpha_{2j+1}$ because $\alpha <
\beta=\alpha_{2n-1}$. Then we define 
$$f_{\alpha}^{\beta}(x)= f_{\alpha}^{\beta}(\alpha_1,\alpha_2,\cdots,
\alpha_{2n-1},   
\alpha_{2n}) = (\alpha_1,\alpha_2,\cdots,
\alpha_{2j},\alpha,\alpha_{2j+2}).$$ 
Clearly,
$(\alpha_1,\alpha_2,\cdots,\alpha_{2j},\alpha,\alpha_{2j+2})$  satisfies
all the three conditions of a point as it is only a subsequence of the
point $x$. Also $(\alpha_1,\alpha_2,\cdots, \alpha_{2j},\alpha,\alpha_{2j+2})$
 has index $\alpha$ hence
 $(\alpha_1,\alpha_2,\cdots,\alpha_{2j},\alpha,\alpha_{2j+2})\in X_{\alpha}$.

\bigskip\noindent Note that in particular, the map $f^{3}_{2}:X_3\to X_2,$
will be as follows.

\bigskip\noindent $(3, \alpha) \mapsto (2,\alpha)\in X_2$ (Case I) \\
$(0, \alpha, 3, \beta)\mapsto (0, \alpha, 2,
\beta) \in X_2$ \\
$(1, \alpha, 3, \beta)\mapsto (1, \alpha, 2, \beta)\in X_2$ \\
$(2, \alpha, 3, \beta)\mapsto (2, \alpha)\in X_2$ (Case I)\\
$(0, \alpha, 1, \beta, 3, \gamma)\mapsto (0, \alpha, 1, \beta, 2, \gamma)\in
X_2$ \\ 
$(0, \alpha, 2, \beta, 3, \gamma)\mapsto (0, \alpha, 2, \beta)\in X_2$\\
$(1, \alpha, 2, \beta, 3, \gamma)\mapsto (1, \alpha, 2, \beta)\in X_2$\\
$(0, \alpha, 1, \beta, 2, \gamma, 3, \delta)\mapsto (0, \alpha, 1, \beta, 2,
\gamma)\in X_2.$\\ All others are as in Case II.

\bigskip\noindent Having defined  $f^{\beta}_{\alpha},\;\;\alpha \leq \beta,$
let us note the following obvious properties of these bonding maps:
\begin{enumerate}
\item[(i).] $f^{\alpha}_{\alpha}:X_{\alpha}\to X_{\alpha}$ is identity.\\
\item[(ii).] Let $f^{\beta}_{\alpha}:X_{\beta}\to X_{\alpha}$ and
  $f_{\beta}^{\gamma}:X_{\gamma}\to
  X_{\beta},\;\;\alpha<\beta<\gamma,\;\;\alpha, \beta, \gamma
  \in\Omega$ be  two bonding maps. Then
  $$f^{\beta}_{\alpha}of_{\beta}^{\gamma}=
  f^{\gamma}_{\alpha}:X_{\gamma}\to X_{\alpha}.$$
\end{enumerate}
Let us verify a specific example of (ii): Let $f^{6}_{4}: X_6\to
X_4$ and $f^{4}_{3}:X_4\to X_3$ 
be the two bonding maps. We want to verify 
$$f^{4}_{3}f^{6}_{4}=f^{6}_{3}.$$ 
For this we choose an arbitrary element of $X_6$
say,\\ $x=(0,\alpha, 1,\beta, 3,\gamma, 5,\delta, 6, \omega )\in X_6,$ 
 where $\alpha>0,\;\; \beta>1,\;\;\gamma>3,\;\;\delta>5,$ and
 $\omega>6.$  Then by 
 the definition of $f^{\beta}_{\alpha}$ we have

\begin{eqnarray*}  
f^{6}_{4}(x) & = & (0,\alpha, 1,\beta,
3,\gamma,5,\delta)=y,\;\;\mbox{say and} \\ 
f^{4}_{3}(y)& = & (0,\alpha, 1,\beta, 3,\gamma). \;\;\; \mbox{Also,} 
\end{eqnarray*}
$$f^{6}_{3}(x) =f^{6}_{3}(0,\alpha, 1,\beta, 3,\gamma,
5,\delta, 6, \omega ))= (0,\alpha, 1,\beta, 3,\gamma).$$ 

\bigskip  Hence $f^{4}_{3}f^{6}_{4}=f^{6}_{3}$.
 
\bigskip\noindent Thus we have the following  inverse system
$\{X_{\alpha},\;\;f_{\alpha}^{\beta};\;\alpha\leq\beta,\;\;\alpha,
\beta\in\Omega\} $ of sets and maps defined on the directed set $\Omega.$

$$\cdots X_{\gamma}\stackrel{f^{\gamma}_{\beta}}\longrightarrow X_{\beta}\stackrel{f_{\alpha}^{\beta}}\longrightarrow X_{\alpha}\cdots X_2\stackrel{f^{2}_{1}}\longrightarrow X_1,$$

\medskip\noindent where $1<2<\cdots <\alpha < \beta < \gamma$ and $1,2,\cdots , \alpha, \beta, \gamma \in\Omega$

\bigskip\noindent Now we verify that the bonding maps
$f^{\beta}_{\alpha}:X_{\beta}\to X_{\alpha}, \;\;\alpha \leq \beta,$
are onto:

\bigskip Let $x=(\alpha_1,\alpha_2,\cdots, \alpha_{2n-1}, \alpha_{2n})
\in X_{\alpha}$ so that $\alpha=\alpha_{2n-1}.$ We choose $\gamma>\beta$ and
then  consider the sequence of even number of elements from the
directed set $\Omega$ and let $y=(\alpha_1,\alpha_2,\cdots,
\alpha_{2n},\beta,\gamma)$. Since $x$ is a subsequence of $y$ and $x$
is a point and also $\alpha_{2n-1}<\alpha<\beta<\gamma$, to
prove that $y$ is a point in  $X_{\beta}$, it suffices to verify that
$\beta\nless\alpha_{2j+1}, \;\; 0\leq j<n$. But if this is true that
$\beta\leq\alpha_{2j+1}, \;\; 0\leq j<n$ will imply that
$\alpha<\alpha_{2j+1}$ since $\alpha < \beta$, which is contradiction to the
fact that $x=(\alpha_1,\alpha_2,\cdots, \alpha_{2n-1},\alpha_{2n})$ is a
point. Thus $y=(\alpha_1,\alpha_2,\cdots,\alpha_{2n},\beta,\gamma)$ is
a point and it is an element of $X_{\beta}$ such that
$f_{\alpha}^{\beta}(y)=x$ by the definition of $f_{\alpha}^{\beta}$. 

\bigskip\noindent Let us  see a particular example of ontoness as follows :

\bigskip Consider $f^{5}_{3}:X_5\to X_3$.  We will show that  there is
 preimage of all the elements 
of $X_3$ (all the  8 types as discussed earlier), is in $X_5$. The preimage of
any element of $X_3,\;\;x=(\alpha_1,\alpha_2,\cdots,\alpha_{2n-1},\alpha_{2n})$
where $\alpha_{2n-1}=3$ can be obtained by just introducing two more
 elements from the directed set $\Omega$ in point $x$ as
$(\alpha_1,\alpha_2,\cdots,\alpha_{2n-1},\alpha_{2n},5,\delta),\;\;\delta>5.$
Thus we see
\begin{eqnarray*}
(f^{5}_{3})^{-1}(x)& = &
  (f^{5}_{3})^{-1}(\alpha_1,\alpha_2,\cdots,\alpha_{2n-1},\alpha_{2n})
  \\
& =&  (\alpha_1,\alpha_2,\cdots,\alpha_{2n-1},\alpha_{2n},5,\delta)\\
& =& y \in X_5
\end{eqnarray*} 
So, we have now an inverse mapping  system based on $\Omega$ 
$$\{X_{\alpha},\;\;f_{\alpha}^{\beta}:X_{\beta}\to
X_{\alpha},\;\;\alpha\leq\beta,\;\;\alpha, \beta\in\Omega\},$$ in which
each $X_{\alpha}\neq \phi$ and each $f_{\alpha}^{\beta}:X_{\beta}\to
X_{\alpha},\;\alpha \leq \beta$ is onto. 

\section{Proof of the main result}

We claim that for the above inverse system $X$, $\varprojlim X_{\alpha}=\phi$.

\bigskip Let us assume the contrary  and suppose there exists an element $x$ in the
inverse limit $\varprojlim X_{\alpha}$. In other words, $x\in \prod X_{\alpha}$
where $$x= (x_0, x_1, \cdots,
x_{\omega},x_{\omega+1},\cdots,x_{2\omega},x_{2\omega+1},
\cdots),\;\;x_i\in X_i \eqno{(*)}$$ such that  whenever $\alpha
\leq\beta,\;\; f_{\alpha}^{\beta}(x_{\beta})=x_{\alpha}$. 

\bigskip Let $\O$ be the set of orders of these $x_i 's$ in $(*)$. It is
clear that this set $\O$ is a cofinal set of $\Omega$, since for any
$\alpha \in \Omega$ there is a set $X_{\alpha}$ and since $x \in \prod
X_{\alpha}$, the i-th co-ordinate in $x$ is some element
$x_{\alpha}\in X_{\alpha}$. Since $\alpha=$index of $x_{\alpha}<$ order of
$x_{\alpha} \in \O$, we find that $\O$ is a cofinal set in $\Omega.$

\bigskip We also observe that  if length of $x_{\alpha}=$ length of
$x_{\beta}$, then order of $x_{\alpha}=$ order of $x_{\beta}$. To prove
this we choose a $\gamma > \alpha$ and $\gamma > \beta$. Then there
exists an element $x_{\gamma}\in X_{\gamma}$ in $(*)$ such that
$f^{\gamma}_{\alpha}(x_{\gamma})=x_{\alpha}$ and
$f^{\gamma}_{\beta}(x_{\gamma})=x_{\beta}$. But from the definition of
the bonding map $f^{\beta}_{\alpha}$ it follows that the orders of
$x_{\alpha}$ and $x_{\beta}$ are some element in the sequence
$x_{\gamma}$, say order of $x_{\alpha}= \alpha_{2i}$ and order of
$x_{\beta}= \alpha_{2j}$. Thus if the length of $x_{\alpha}= i=$
length of $x_{\beta}= j,$  then clearly, $\alpha_{2i}=\alpha_{2j}$
i.e., order of $x_{\alpha}=$ order of $x_{\beta}$.

\bigskip Therefore, the  set of orders $\O$ behaves according to the lengths
of $x_i$ in $(*)$, and length of any point is a natural number.

\bigskip Thus, if the lengths of $x_i 's$ are unbounded then we will
get a simple sequence of the orders of the elements $x_i$ which will
be cofinal. Hence there will exist a cofinal simple 
sequence in $\Omega=[0,\omega_1)$. But this is clearly a contradiction.
 On the other hand if the lengths of
  $x_i's$ are bounded, then the set of orders of $x_i's, \;\;i.e., \;\;\O$
  will contain a maximal element of $\Omega=[0,\omega_1)$ which is again a contradiction. Hence, we
    find that in either of the two cases viz., when the lengths 
are unbounded or bounded we have a contradiction since it is well
known that  the set
$\Omega=[0,\omega_1)$ neither posseses a simple cofinal sequence nor a
  maximal element. Hence there can not exist any element in the
  inverse limit, i.e., $\varprojlim X_{\alpha}=\phi$.\hfill\rule{2mm}{2mm}

{\Remark \rm In view of the above construction, it is clear that one can always have an inverse system in any category (e.g., topological spaces and continuous functions or modules and module homomorphisms etc.) admitting arbitrary products with nonempty objects and onto bonding morphism whose inverse limit can be empty.}

\bigskip \noindent
Satya Deo, \\
Harish Chandra Research Institute, \\
Chhatnag Road, Jhusi,\\
Allahabad 211 019, India. \\
Email: sdeo@mri.ernet.in,  vcsdeo@yahoo.com

\medskip \noindent
Veerendra Vikram Awasthi, \\
Institute of Mathematical Sciences, \\
CIT Campus, Taramani \\
Chennai 600 013, India. \\
Email: vvawasthi@imsc.res.in, vvawasthi@yahoo.com

\end{document}